\documentclass{amsart}
\usepackage{amssymb}
\usepackage{pb-diagram}

\newcommand\datver[1]{\def\datverp%
{\par\boxed{\boxed{\text{Version: #1; Run: \today}}}}}
\datver{3.0A; Revised: 9/23/1998}

\newcommand\C{\mathbb C}

\newcommand\R{\mathbb R}
\renewcommand\P{\mathbb P}
\newcommand\Z{\mathbb Z}
\newcommand\T{\mathbb T}
\newcommand\ts{{\otimes}}
\newcommand\si{\sigma}

\newcommand\cohom{\operatorname{H}}

\newcommand\De{{\Delta}}

\newcommand\Tn{\mathbb T_{p}}
\newcommand\pa{\partial}

\newcommand\maL{\ell}
\newcommand\mal{\ell}

\newcommand\Indo{\operatorname{Ind}_{V_\Omega}}
\newcommand\Ind{\operatorname{Ind}}
\newcommand\im{\operatorname{Im}}
\newcommand\Id{\operatorname{Id}}
\newcommand\Tr{\operatorname{Tr}}
\newcommand\ar{\rightarrow}
\newcommand\inj{\hookrightarrow}

\newcommand\zp{{\Z^p}}
\newcommand\zo{{\Z[\mu]}}
\newcommand\rp{{\R^p}}
\renewcommand\O{{\Omega}}
\newcommand\ro{{\R^p\times\Omega}}
\renewcommand\o{{\omega}}
\newcommand\tn{{\tau^\mu}}

\newcommand\Co{{C(\Omega)}}
\newcommand\Cz{{C(\Omega,\Z)}}

\newcommand\Cp{{C(\Omega)\rtimes\Z^p}}
\newcommand\V{{V_\Omega}}
\newcommand\Cv{{C(V_\Omega)}}

\newcommand\Mos{M_n(C^{\infty,0}(V_\Omega))}
\newcommand\MC{M_n(\C)}

\newcommand\coinv{{C(\Omega,\Z)_{\zp}}}
\newcommand\coinvd{{C(\Omega,\Z)_{\Z^{2}}}}
\newcommand\tk{{\tau^\mu_*}}
\newcommand\Kp{K_0(C(\Omega)\rtimes\Z^p)}
\newcommand\Kd{K_0(C(\Omega)\rtimes\Z^2)}
\newcommand\Kv{K_0(C(V_\Omega))}
\newcommand\Kvs{K_0(C^{\infty,0}(V_\Omega))}

\newcommand\mt{{\mu_{\O}^\zp}}
\newcommand\chl{{\mbox{ch}_{l}}}

 \newcommand\di{{\partial}}
 \newcommand\de{{\partial_{\O,\rp}^e}}
 \newcommand\I{{\mbox{Ind}_{\V}\partial_{\O,\rp}^e}}

\newcommand\lc{\Omega^*_{\maL,c}}

\newcommand\lle{\Omega^e_{\maL}(\V,\R)}
\newcommand\llk{\Omega^k_{\maL}(\V,\R)}
\newcommand\llp{\Omega^p_{\maL}(\V,\R)}

\newcommand\hlk{\cohom^k_{\maL}(\V,\R)}
\newcommand\hlp{\cohom^p_{\maL}(\V,\R)}
\newcommand\hl{\cohom^*_{\maL}(\V,\R)}

\newcommand\hc{\cohom^*(\V,\R)}

 \newcommand\hle{\cohom^e_{\maL}(\V,\R)}
\newcommand\hlo{\cohom^o_{\maL}(\V,\R)}

\newcommand\luk{\Omega^k_{\maL,c}(U,\R)}
\newcommand\lu{\Omega^*_{\maL,c}(U,\R)}
\newcommand\hulk{\cohom^k_{\maL,c}(U,\R)}
\newcommand\hul{\cohom^*_{\maL,c}(U,\R)}

\newcommand\hule{\cohom^e_{\maL,c}(U,\R)}
\newcommand\hulo{\cohom^o_{\maL,c}(U,\R)}

\newcommand\hwl{\cohom^*_{\maL,c}(V,\R)}

\newcommand\luu{\Omega^*_{\maL,c}(U_0,U_1,\R)}
\newcommand\huulk{\cohom^k_{\maL,c}(U_0,U_1,\R)}
\newcommand\huul{\cohom^*_{\maL,c}(U_0,U_1,\R)}

\newtheorem{theorem}{Theorem}

\newtheorem{proposition}{Proposition}
\newtheorem{corollary}{Corollary}
\newtheorem{conj}{Conjecture}

\theoremstyle{definition}
\newtheorem{definition}{Definition}

\theoremstyle{remark}
\newtheorem{remark}{Remark}

\begin{document}

\title[Gap-labelling for  quasi-crystals]{Gap-labelling for
quasi-crystals\\ (Proving  a conjecture by J. Bellissard)}
\author[M. Benameur]{Moulay-Tahar Benameur}
\address{Inst. Desargues, Lyon, France}
\email{benameur@desargues.univ-lyon1.fr}
\author[H. Oyono]{Herv\'e Oyono-Oyono}
\address{Universit\'e de Clermont-Ferrand}
\email{herve.oyono@math.univ-bpclermont.fr}
\begin{abstract}
In the present paper, we summarize  a proof of the Bellissard gap-labelling conjecture for quasi-crystals. Our main tools are the measured index theorem for laminations together with the naturality of the longitudinal Chern character.
\end{abstract}
\maketitle

\vspace*{1cm}
\section{Introduction}

Let us consider a Schrodinger operator on $\R^p\,$
$$
 H=-\frac{h^2}{2m}\De+V\,,
$$
where $\De$ is the Laplacian on $\R^p$ and $V$ is a potential.
The set of observables affiliated to this Schrodinger operator is a
$C^*-$ algebra. It must contain at least the $C^*-$algebra generated
by operators $(H-z\Id)^{-1}$ where $z$
belongs to the resolvent of $H\,$.
If $H$ describes a particle in a homogeneous media, the
physical properties of this media do not depend upon
the choice of an origin in $\R^p\,$. In particular, the algebra of
observables must also contain the $C^*$-algebra $C^*(H)$ generated by
$T_a(H-z\Id)^{-1}T_{-a}$ where $T_a$ is
the operator of translation by  $a\in\rp\,$.

In \cite{Bellissard}, J. Bellissard has attached to this Schrodinger
operator $H$ a compact space $\O_H$  equipped with a minimal  action
of
$\rp$
such that the crossed product algebra $C(\O_H)\rtimes\rp$
 contains $C^*(H)\,$.
The space $\O_H$ is called the hull of $H$ and  is by definition the
strong
closure in $\mathcal{B}(L^2(\rp))$
of the family
$$
 \{T_a(H-z\Id)^{-1}T_{-a};\, a\in\rp\}\,
$$
where   we have fixed an element  $z$  in  the
resolvent of $H\,$. The action of $\rp$ on $\O_H$ is induced by
translations and up to a $\rp-$equivariant homeomorphism, the space
$\O_H$ is independant on the choice of $z$ in the resolvent of $H\,$.


A typical potential for motion of conduction electrons is given by:
$$
 V(x)=\sum_{y\in L}v(x-y)\,,
$$
where
$L$ is the point set of equilibrium positions of atoms and $v$ is the
effective potential for a valence electron near an atom (see
\cite{BellissardHerrmannZarrouati}). This set $L$ is {\it uniformly
discrete}, i.e.
there exists a positive number  $r$ such that any ball of radius $r$
contains at most one
point.
In \cite{BellissardHerrmannZarrouati}, J. Bellissard, D. J. L.
Hermann and M. Zarrouati had attached to this point set
$L$ a compact space $\O_L\,$, called the hull of $L$ equipped with  a
minimal
action of $\rp\,$. The space  $\O_L$ is called the hull of $L$
and is constructed as follow: Let $\nu_L$ be the measure on $\rp$
defined
for every compactly supported continuous
fonction $f\,$ by $\nu_L(f)=\sum_{y\in L}f(y)\,$. Then $\O_L$ is by
defininition the
weak-$*$ closure (with respect to the compactly supported continuous
functions on $\rp$) of the family of translations of $\nu_L$ by the
elements
of $\rp\,$.  Note that  $V=\nu_L*v$ and more generally,  for
$\nu\in\O_L$ we set:
$$
 H_\nu=-\frac{h^2}{2m}\De+\nu*v\,.
$$
 For all $z$ in the resolvent of the operator $H\,$,  the
map
$$
 \nu\in\O_L\mapsto (H_\nu-z\Id)^{-1}\in\O_H
$$
 is continuous, equivariant  and surjective
(see \cite{BellissardHerrmannZarrouati}). Thus, the crossed product
algebra
$C(\O_H)\rtimes\rp$ lies in the crossed product algebra
$C(\O_L)\rtimes\rp$ and in particular $C(\O_L)\rtimes\rp$ contains
the algebra $C^*(H)\,$.
The main advantage of dealing with $\O_L$ rather than $\O_H$ is that
$\O_L$ only depends on the geometry of $L\,$; For instance, if
$L$ is given by a rank $p$ lattice $\mathcal{R}$ in $\rp\,$, then the hull of
$L$
is  $\rp/\mathcal{R}\,$.  Actually in this
case $C^*(H)$ can be computed by using Bloch theory \cite{Bellissard}: we can check
that $C^*(H)=C(B)\ts \mathcal{K}\,$, where
$\mathcal{K}$
is the algebra of compact operator on a separable Hilbert space and
$B$ is the Brillouin zone, defined by $B=\rp/\mathcal{R}^*\,$,
where $\mathcal{R}^*$ is the reciprocal lattice of $\mathcal{R}\,$.


To the Schrodinger operator $H$ is associated the {\it  integrated
density of states} (see \cite{Bellissard}) $E\mapsto \mathcal{N}(E)\,$,
where $\mathcal{N}(E)$   is defined as the number of states by
unit of volume with egenvalues less or equal to $E\,$.
The remarkable result of J. Bellissard, D. J. L.
Hermann and M. Zarrouati in  \cite{BellissardHerrmannZarrouati} is
that  the   integrated
density of states on gaps must take value in a countable subgroup
of $\R$ that only depends  on the point set $L\,$. This computation
goes as follows.
If $\P$ is a $\rp-$invariant probability measure on $\O_L\,$, then
$\P$
induces a trace
$\tau^\P$ on $C(\O_L)\rtimes\rp$ and this trace extends to the
von-Neumann algebra of $C(\O_L)\rtimes\rp\,$. Let us denote by
$\chi_{]-\infty,E]}$ the characteristic function of the set
$]-\infty,E]\,$. Then $\chi_{]-\infty,E]}(H)$ belongs to the
von-Neumann
algebra of $C(\O_L)\rtimes\rp\,$.
For an ergotic probability $\P\,$, J. Bellissard stated in
\cite{Bellissard}   the so-called
{\it Shubin's formula} for $H\,$:
$$
        \mathcal{N}(E)=\tau^\P(\chi_{]-\infty,E]}(H))\,.
$$

In particular, if $E$ is in a spectral gap of $H\,$, then
$\chi_{]-\infty,E]}$
is a continuous function on the spectrum of $H$ and thus,
$\chi_{]-\infty,E]}(H)$ is an idempotent  in $C(\O_L)\rtimes\rp\,$
(recall
that
$H$ is bounded below). In consequence, according  to  the Shubin
formula,
the value of  $\mathcal{N}$ on spectral gaps of $H$ belongs to the
image
of the additive map
$$
 \tau^\P_*:K_0(C(\O_L)\rtimes\rp)\to\R\,,
$$
where $\tau^\P_*$ is the
morphism induced by the trace $\tau^\P$ in $K-$theory.


We now focus on the case where the point set  $L$ is   a quasicrystal which is obtain by {\it
the
cut-and -project method}
(see \cite{DuneauKatz} and also
\cite{BellissardHerrmannZarrouati}).
To such a point set, J. Bellissard, E. Contensous and A. Legrand
have associated in \cite{BellissardContensousLegrand} (see also
\cite{BellissardHerrmannZarrouati})  a  canonical minimal dynamical
system
$(\T_L,\zp)\,$
which is Morita equivalent to  $(\O_L,\rp)$ and such that $\T_L$  is a
Cantor set. Moreover, there is a canonical
ergodic invariant probability measure $\mu$  on $\T_L$ such that
$$
 \tau^\P_*(K_0(C(\O_L)\rtimes\rp))=\tau^\mu_*(K_0(C(\T_L)\rtimes\zp)),
$$
where $\tau^\mu$ is the trace on $C(\T_L)\rtimes\zp$ induced by
$\mu\,$. Eventually, the image
$$\tau^\P_*(K_0(C(\O_L)\rtimes\rp))$$ is predicted by the
following  conjecture of J. Bellissard \cite{BellissardKellendonkLegrand}:

Let $\O$ be a Cantor set equipped with an
action of
$\zp$ and with a $\zp$-invariant
 measure $\mu\,$.
The measure $\mu$ induces
a trace $\tn$ on the crossed product $C^*$-algebra   $\Cp\,$.
Let us denote   by $\zo$ the additive subgroup of $\R$ generated
by  $\mu$-measures of compact-open subsets of $\O\,$.
We make the assumption that $\O$ has no non-trivial compact-open
invariant
subsets (this is clearly the case if the action of $\zp$ is minimal).
\begin{conj}\label{conjecture}
{\bf (The Bellissard gap-labelling conjecture)}

$$
        \tk(\Kp)=\zo.
$$
\end{conj}
The goal of the present paper is to give a proof of
the Bellissard conjecture. The method adopted uses the measured index
theorem for  laminations as proved in  \cite{MooreSchochet}.

{\em Acknowledgements.} We are indebted to J. Bellissard for
explaining the gap-labelling problem  to us  and
for his constant encouragements. We
also would like to thank E. Contensous, T. Fack, J. Kellendonk, A. Legrand,
V. Nistor and  C. Schochet   for
several helpful discussions.

\section{The mapping torus}
Let $\O$ be a Cantor set. Assume that the group
$\zp$ acts on $\O$ by homeomorphisms and that there exists a
$\zp$-invariant
 measure $\mu$ on $\O$.  We
assume that $\O$ has no compact-open $\zp$-invariant subset except
$\emptyset$ and $\O\,$. This is the case in particular if the action
of $\zp$ on $\O$ is minimal.

The action of $\zp$ on $\O$ induces an
action of $\zp$ on the $C^*$-algebra $\Co$  and thus, we can form the
crossed product  $C^*$-algebra $\Cp$. The measure $\mu$ induces
a trace $\tn$ on the $C^*$-algebra   $\Cp$
and we obtain in this way a group morphism
       $ \tk:K_0(\Cp) \to \R\,$.
In what follows, we shall denote by $\zo$ the additive subgroup of $\R$
generated
by the $\mu$-measures of compact-open subsets of $\O$
and by $\coinv$ the coinvariants of the action of
$\Z^p$ on $\Cz\,$, i.e. the
quotient of $\Cz$ by the subgroup generated by elements
of the form $n(f)-f$, where $f\in\Cz$ and $n\in\zp$.

We are interested in computing the image  of $K_0(\Cp)$ under
$\tau^{\mu}_*$.

In the case $p=1\,$,
E. Contensous  proved in \cite{Contensous}, using the
Pimsner-Voiculescu six term exact sequence,  that the inclusion
$\Co\inj C(\Omega)\rtimes\Z$ induces an isomorphism
$C(\Omega,\Z)_{\Z}\simeq K_0(C(\Omega)\rtimes\Z)$ and thus,
$$\tk(K_0(C(\Omega)\rtimes\Z)) = \zo\,.$$

In the case $p=2\,$, the computation of $\Kd$ was carried out by
J. Bellissard, E. Contensous and A. Legrand in \cite{BellissardContensousLegrand}
using the Kasparov spectral
sequence.
The result is:
$$
        \Kd\simeq \coinvd \oplus \Z.
$$
More precisely:
\begin{itemize}
        \item   The inclusion
                $\Z\inj K_0(C(\Omega)\rtimes\Z^2)$ maps the canonical
generator of
                $\Z$ to the image, under the morphism induced by the
inclusion
                $C^{*}_{r}(\Z^2)\ar C(\Omega)\rtimes\Z^2\,$, of the
Bott generator in the
                $K$-theory of the $C^*$-algebra $C^{*}_{r}(\Z^2)\simeq
                C(\T^{2})$ of the discrete group $\Z^2$.
        \item   The inclusion $\coinvd\inj K_0(C(\Omega)\rtimes\Z^2)$
is induced by
                the inclusion $\Co\inj C(\Omega)\rtimes\Z^2$.
\end{itemize}

Recall that the Bott generator of  $K_0(C(\T^{2}))$ is the unique
element with Chern character egal to the volume form of $\T^2\,$. In
particular, it is traceless and therefore,
the above computation gives:
$$
        \tk(K_0(C(\Omega)\rtimes\Z^2)) = \zo.
$$

In the case $p=3$, the computation of the image  of $K_0(C(\O)\rtimes
\Z^3)$ under the
trace has been recently
performed by J. Bellissard, J. Kellendonk and A. Legrand
\cite{BellissardKellendonkLegrand} and gave the same result.
These computations
lead J. Bellissard to state  conjecture \ref{conjecture}. We can 
without loss of generality assume that $p$ is even (see 
\cite{BenameurOyono}). From now on,
$p$ will denote an even integer.

The group $K_0(C(\Omega)\rtimes\Z^p)$ can be computed in
term of Kasparov cycles out of the mapping torus isomorphism.
We recall in the end of this section the construction of this isomorphism.

The mapping torus is by definition the space
$$\V=\frac{\O\times\rp}{\zp}\,.$$
This space is a foliated space with leaves given by the projections of
$\{\o\}\times\rp$ in $\O\times\rp\to\V\,$, where
$\o$ is an element of $\O\,$.

Let
$\di : C^\infty(\R^p)\ts S^+\to  C^\infty(\R^p)\ts S^-$ be the Dirac
operator on $\rp$
where $S^+$ and $S^-$ are the two irreducible spin representations of
$\mbox{Spin}(p)\,$. The Dirac operator $\di$ is $\zp-$equivariant.
Let $e$ be a projection of $M_n(C(\V))$ and let $\tilde{e}$ be the
$\zp-$invariant projection of $C(\O\times\rp)\otimes M_n(\C)$
corresponding to $e$ under the  projection $\O\times\rp\to\V\,$.
If we assume $\tilde{e}$ smooth in the $\rp-$direction, then
the operator
$$
        \tilde{e}\circ (Id_{C(\O)}\otimes \di\otimes I_n): \tilde{e}(C(\O)\ts
        C^\infty(\R^p)\ts S^+\otimes\C^n)\to  \tilde{e}(C(\O)\ts
        C^\infty(\R^p)\ts S^-\otimes\C^n)
$$
is a $\zp-$equivariant elliptic differential operator in the
$\rp-$direction.
Hence it induces a longitudinal    elliptic       differential operator
$\de$
on $\V\,$. According to \cite{MooreSchochet}, the operator $\de$
admits a $K-$theory index $\I$  which belongs to
the $K$-theory group
$\Kp\,$.

We can easily  check that the map $e \mapsto \Indo(\de)$ induces a well
defined morphism from
 $K_0(C_0(\V))$ to $K_0(C(\Omega)\rtimes\Z^p)\,$.
The following theorem is due to A. Connes (see
\cite{ConnesBook})

\begin{theorem}\label{torus}\
 The  map   $e\mapsto\I$ induces an isomorphism:
$$
        \mt: \Kv {\stackrel{\simeq}{\longrightarrow}} \Kp.
$$
\end{theorem}

 In consequence of this
 theorem, it is enought for proving the Bellissard conjecture to
 check that $\tau^\mu_*(\I)$ belongs to $\zo\,$.

 Assume first that $\O=\{\mbox{pt}\}$ is just a point. The mapping
torus is then
 $\V=\T_p\,$, the usual $p-$torus. If $e$ is a smooth projector in
$M_n(\T_p)\,$, then
 $\partial_{\{\mbox{pt}\},\rp}^e$ is just the Dirac operator of the
$p-$torus $\T_p$ with
 coefficients in the vector bundle over $\T_p$ associated with $e\,$.
 The Atiyah index formula for coverings (see \cite{Atiyah}) asserts
 then that
 $$
        \tau_*(\mbox{Ind}\,\, \partial_{\{\mbox{pt}\},\rp}^e)= \langle
        \mbox{ch}([e]),\,[\T_p]\rangle \,
$$
 where $\tau$ is the canonical trace on $C^*(\Z^p)\,$,  $[\T_p]\in H_p(\T_p,\R)$ is the fundamental class
of $\T_p$ and
 $\mbox{ch}([e])\in H^*(\T_p,\R)$ is the Chern character of
 $[e]\in \Kv\,$.

 The Connes measured longitudinal index theorem  (see
 \cite{ConnesSurvey}), or rather its extension to foliated spaces
 (see \cite{MooreSchochet}) provides such a formula for $\de\,$.
 To state this formula for $\V\,$, we need to define suitable
 characteristic classes and a cycle to integrate these classes.
  This will be done in the two
following sections.
 The proof of the Bellissard conjecture will then be given in Section
 6.

 \section{Longitudinal characteristic classes}
 The support for the characteristic classes involved in  the measure
index theorem for foliation
 is the longitudinal cohomology of $\V\,$. We give in this section,
 first the definition of this cohomology and then, the construction of
the
 longitudinal Chern character on $\Kv$ valued in this longitudinal cohomology.

 \subsection{The longitudinal de Rham complex}
 Let $\O^k(\R^p)$ be the space of $k$-differential forms on the vector
space $\rp\,$, endowed  with its usual Frechet
topology.

\begin{definition}
A (real) longitudinal differential $k$-form on $\V$ is a
$\zp$-equivariant
continuous map $\phi:\O\ar \O^k(\R^p)$.
\end{definition}

We denote by  $\llk\,$ the space of longitudinal differential $k-$forms
on $\V\,$.
If $\phi$ is a longitudinal $k-$form, its longitudinal differential
$d_{\maL}(\phi)$ is by
definition the
longitudinal $(k+1)$-differential form which is given by the
$\Z^p$-equivariant map $\omega \mapsto
d(\phi(\omega))$ where $d$ is the de Rham differential
on  $\R^p$. 
Hence:
$$
        d_{\maL}:\llk\ar\O^{k+1}_{\maL}(\V,\R)
$$
provides a differential structure on the graded vector space
$\O^*_{\maL}(\V,\R)=\bigoplus\llk\,$
and satisfies $d_{\maL}\circ d_{\maL}=0$. In what follows the
superscripts $^e$ and $^o$ mean
respectively
even and odd forms or classes of forms.
The cohomology of the complex $(\O^*_{\maL}(\V,\R), d_{\maL})$ will be
denoted by
$$
        \hl=\bigoplus_{k\geq 0} \hlk = \hle \bigoplus \hlo.
$$

\begin{remark}\label{cech} We can check (see  \cite{MooreSchochet})
that $\hl$ is
  the cohomology of the sheaf of continuous functions which are locally
constant
  in the leaf direction or equivalently the sheaf of continuous and
  equivariant functions on equivariant open subsets of
  $\O\times\R^p\,$, constant in the $\R^p-$direction.
If $\hc$ denotes the Cech cohomology groups of $\V$ with
real  coefficients, then we have a well defined morphism:
$$
        \hc \ar \hl
$$
induced by the natural morphism of sheaves.
\end{remark}

We define the support of a longitudinal form  $\phi$ as the image of
the support of $(\o, t)\mapsto
\phi(\o, t)$ under the projection $\O\times\rp\to\V\,$. The support of
$\phi$ is thus a compact subset of $\V\,$. Let now $U$ be an open
subset  of $\V\,$ and let
$\lu=\bigoplus\luk\,$ be the space of longitudinal
differential $k-$forms with support in $U$.

The restriction of the longitudinal differential $d_{\maL}$ to $\lu$
preserves it and
 $(\O^*_{\maL,c}(U,\R), d_{\maL})$ is a subcomplex  of the longitudinal
complex
$(\O^*_{\maL}(\V,\R), d_{\maL})$.
We denote  by
$$
        \hul=\bigoplus\hulk = \hule \bigoplus \hulo
$$
the cohomology of this subcomplex.
Similarily, we can also define  the relative longitudinal cohomology
for
relative open pairs.
A relative open pair $(U_0,U_1)$ is by definition given by two open
subsets $U_0$ and $U_1$ of $\V$
such that $U_0\subset U_1\,$. From such a pair $(U_0,U_1)\,$, we
can define the differential complex $(\luu, d_{\maL})$ as the
quotient differential complex in the exact sequence:

$$
        0 \ar \, \lc(U_0,\R) \ar \, \lc(U_1,\R) \ar \, \luu \ar 0.
$$
The cohomology of the resulting complex
$(\luu\,, d_{\maL})$
is called the relative longitudinal cohomology of the relative open
pair $(U_0,U_1)$ and is denoted by:
$$
        \huul=\bigoplus\huulk = \cohom^e_{\maL,c}(U_0,U_1,\R)
\bigoplus
\cohom^o_{\maL,c}(U_0,U_1,\R).
$$

The following proposition is a classical  homological algebra result,
see for
instance \cite{Spanier}.

\begin{proposition}\label{longexactsequence} The following long exact
sequence in longitudinal
cohomology holds:
$$
        \cdots\ar \cohom^k_{\maL,c}(U_0,\R)\ar
\cohom^k_{\maL,c}(U_1,\R)\ar \huulk\ar
        \cohom^{k+1}_{\maL,c}(U_0,\R)\ar\cdots
$$

\end{proposition}

\subsection{The longitudinal Chern character}

We give in this subsection the  definition and main properties of the
longitudinal Chern character:
$$
        ch_{\maL} : K^*(\V) \to \cohom_{\maL}^{[*]}(\V,\R), \quad
*=0,1,
$$
where $\cohom_{\maL}^{[*]} = \oplus_{j\in \Z}
\cohom_{\maL}^{*+2j}(\V,\R)$.

\begin{definition} A continuous function $f$ on $\V$ is longitudinally
smooth if, viewed as a
$\zp-$invariant
function on $\ro\,$, it is smooth in the $\rp-$direction (i.e an
element of $\O^\maL_0(\V,\R)\otimes\C\,$). We denote by
$C^{\infty,0}(\V)$ the algebra of continuous longitudinally smooth
functions
on $\V\,$.
\end{definition}
The algebra $C^{\infty,0}(\V)$ is a dense subalgebra
of the algebra $\Cv$ of
continuous functions on $\V$,
which is stable under holomorphic functional calculus. Hence the
inclusion
$i: C^{\infty,0}(\V) \hookrightarrow \Cv$ induces an isomorphism
\cite{ConnesBook}:
$$
        i_*: \Kvs \stackrel{\simeq}{\rightarrow}\Kv\,.
$$

    Let $e$ be a projector in $\Mos\,$. Let us denote by ${\tilde e}$
    the smooth in the $\rp-$direction $\zp-$invariant
$\MC\,$-valued map on $\ro\,$,
defined by the projection $e\,$.
Then $\displaystyle \Tr\,{\tilde e}\exp\left(\frac{d_{\maL}{\tilde e}d_{\maL}{\tilde
e}}{2i\pi}\right)$  is an element of $\lle\,$.
We have the following proposition:
\begin{proposition}\cite{BenameurOyono}
If $e$ is a projector in $\Mos\,$, then
\begin{enumerate}
\item  $\displaystyle \Tr\,{\tilde e}\exp\left(\frac{d_{\maL}{\tilde e}d_{\maL}{\tilde
e}}{2i\pi}\right)$ is a $d_{\maL}$-closed differential form;
\item  If $[e]$ and $[e']$ define the same class in $\Kvs\,$, then
$$\Tr\,{\tilde e}\exp\left(\frac{d_{\maL}{\tilde e}d_{\maL}{\tilde
e}}{2i\pi}\right)-\Tr\,{\tilde e'}\exp\left(\frac{d_{\maL}{\tilde e'}d_{\maL}{\tilde
e'}}{2i\pi}\right)$$ is a coboundary of $\lle\,$.
\end{enumerate}
\end{proposition}
As a consequence of this proposition, the class of
 $\displaystyle\Tr\,\tilde e\exp\left(\frac{d_{l}{\tilde
 e}d_{l}{\tilde e}}{2i\pi}\right)$ in $\hl$  only
 depends on the class of $e$ in $\Kvs$. 

 \begin{corollary}
     There is a morphism $\chl:\Kv\ar \cohom_{\maL}^{e}(\V,\R)$ such
that for every
     smooth projector $e$ of $\Mos\,$,
     $\displaystyle \chl([e])=\Tr\,
     {\tilde e}\exp\left(\frac{d_{l}{\tilde e}d_{l}{\tilde
e}}{2i\pi}\right)\,$.
     \end{corollary}
      The morphism $\chl$ is called {\it the longitudinal chern
     character}.
     The longitudinal Chern character can  also be defined for the odd
 $K-$theory of $\Cv\,$. Hence we obtain in this way a morphism
\begin{multline*}
        \chl:K_*(\Cv)=\Kv\oplus K_1(\Cv)\longrightarrow \\
\cohom_{\maL}^{*}(\V,\R)=\cohom_{\maL}^{e}(\V,\R)\oplus\cohom_{\maL}^{o}(\V,\R)\,.
\end{multline*}
      If $U$ is an open subset of $\V\,$, we can in the same way
      define  a longitudinal Chern character $\chl:K_*(C_0(U))\ar
      \cohom_{\maL,c}^{*}(U,\R)\,$. Moreover, the longitudinal Chern
character
      admits a relative version: for every relative open pair
$(U_0,U_1)$
      of $\V\,$, there is a map $\chl:K_*(C_0(U_0\setminus U_1))\ar
      \cohom_{\maL,c}^{*}(U_0,U_1,\R)\,$.
      The crucial point for the proof of the Bellissard conjecture is
      that this Chern character is compatible with the long exact
sequences
      associated with  a relative open pair $(U_0,U_1)\,$ of $\V$
\cite{BenameurOyono}.

      \begin{theorem}\label{intertwin}
The longitudinal Chern character intertwines the two following  exact
sequences:
$$
        \begin{diagram}
        \node{K_{0}(C_0^{\infty,0}(U_{0}))}\arrow{e}
        \node{K_{0}(C_0^{\infty,0}(U_{1}))}\arrow{e}
        \node{K_{0}(C_0^{\infty,0}(U_{0},U_{1}))}\arrow{s}\\
        \node{K_{1}(C_0^{\infty,0}(U_{0},U_{1}))}\arrow{n}
        \node{K_{1}(C_0^{\infty,0}(U_{1}))}\arrow{w}
        \node{K_{1}(C_0^{\infty,0}(U_{0}))}\arrow{w}
        \end{diagram}
$$
and
$$
        \begin{diagram}
        \node{ \cohom^e_{l,c}(U_{0},\R)}\arrow{e}
        \node{ \cohom^e_{l,c}(U_{1},\R)}\arrow{e}
        \node{ \cohom^e_{l,c}(U_{0},U_{1},\R)}\arrow{s}\\
        \node{ \cohom^o_{l,c}(U_{0},U_{1},\R)}\arrow{n}
        \node{ \cohom^o_{l,c}(U_{1},\R)}\arrow{w}
        \node{ \cohom^o_{l,c}(U_{0},\R))}\arrow{w}
        \end{diagram}
$$
\end{theorem}

\section{The fundamental cycle}

We give in this section the definition of the Ruelle-Sullivan cycle
associated with the
$\Z^p$-invariant measure $\mu$ on the foliated bundle $\V\,$. This
cycle
allows to integrate the longitudinal top dimensional classes.
 Let $\chi$ be  the characteristic function
 of the open  set $U=]0,1[^p$ in $\R^p$, we define $C_{\mu,\Z^p}:\llp
\to \R$, by:
$$
       C_{\mu,\Z^p}(\phi) :=
        \left<\mu\otimes\int_{\R^p},\,  \chi\phi \right>.
$$
The next proposition shows that  actually, the map $C_{\mu,\Z^p}$
induces a well
defined
map on $\hlp\;$  \cite{BenameurOyono}:

\begin{proposition}\
If $\phi$ is a $d_\ell-$coboundary, then
$<C_{\mu,\Z^p}, \phi>=0$.
\end{proposition}

We are now in position to state the measured index theorem for the
longitudinal Dirac operator with coefficients in a continuous
longitudinally smooth vector bundle \cite{BenameurOyono}.

\begin{theorem}\label{measured.index}\ The measured index of the
longitudinal Dirac
operator $\de$
with coefficients in the vector bundle $E$ associated with a
projector $e$ of $M_n(\Cv)$ is given by:
$$
        \Ind_{\Z^p,\mu}(\de) =  <ch_{\maL}([e]), [C_{\mu, \Z^p}]>.
$$
\end{theorem}

We end this section by identifying  the top dimension 
group of longitudinal cohomology with the coinvariants of the algebra 
$C(\O,\R)\,$ of real valued continuous fonction on $\O\,$.
Recall that   the coinvariants $C(\O,\R)_\zp$ of the action of
$\Z^p$ on   $C(\O,\R)$ is the
quotient of $C(\O,\R)$ by the subgroup generated by elements
of the form $n(f)-f$, where $f\in C(\O,\R)$ and $n\in\zp$. 
\begin{proposition}\ \cite{BenameurOyono}
If $\phi$ is an element of $\O^p_{\maL}(\V,\R)$, we define the continuous function
$\Psi_{\Z^p}(\phi)$ of $C(\O,\R)$ by setting:
$$
        \Psi_{\Z^p}(\phi)(\o) := \int_{]0,1[^p} \phi(\o, {x}) d{x}, \quad
\forall \o\in \O.
$$
Then,  $\phi\to \Psi_{\Z^p}(\phi)$ induces a map from $\hlp$ to the 
coinvariants $C(\O,\R)_\zp\,$.
\end{proposition}

Moreover, this map identifies $\cohom^{p}_{\maL}(\V, \R)$ with
$C(\O,\R)_\zp\,$:

\begin{theorem}\label{isocoinv} \cite{BenameurOyono}
The transform $\psi_{\Z^p}$ is an isomorphism, i.e.
$$
        \cohom^{p}_{\maL}(\V, \R)  \cong C(\O,\R)_\zp .
$$
\end{theorem}

Since the 
measure $\mu$ is $\zp-$invariant, it induces a linear form on $C(\O,\R)_\zp$ that 
we shall denote by $\Phi_{\Z^p, \mu}\,$. Under the above 
identification between $C(\O,\R)_\zp\,$ and $\cohom^{p}_{\maL}(\V, \R)\,$, 
the map 
$\Phi_{\Z^p, \mu}\,$ corresponds
to the Ruelle-Sullivan cycle $C_{\mu,\Z^p}\,$:

\begin{proposition}\ \cite{BenameurOyono}
The following diagram is commutative:
$$
        \begin{array}{ccc}
        \hlp & \stackrel{\Psi_{\Z^p}}{\cong} &C(\O,\R)_\zp\\
         C_{\mu,\Z^p} \searrow & & \swarrow \Phi_{\Z^p, \mu}\\
         & \R &
         \end{array}
$$
\end{proposition}

\section{Proof of the Bellissard conjecture}
We give in this section a proof of the Bellissard conjecture.
We first recall the construction of the Kasparov spectral sequence
associated with the mapping torus $\V$ (see \cite{KasparovInv}).

Let again $\Tn$ be the $p-$torus and let us denote for $j\geq 1$ by $D_j$ the
unit open disk in $\R^j$.  Let $Y_{0}\subset Y_{1}\subset\cdots\subset Y_{p}=\Tn\,$
be a filtration of $\Tn$ by the
skeletons of some triangulation of $\Tn$, where $Y_i$ is the
$i-$skeleton of $\Tn\,$.   We associate to $(Y_{i})_{0\leq i \leq p}$ an open filtration
$(Z_{i})_{0\leq i \leq p}$ of $\Tn$ in the following way:

\begin{itemize}
        \item $Z_{0}$ is an open neighbourhood of $Y_{0}$ which is homeomorphic to
        $Y_{0}\times D_{p}$.
        \item Assume that for $1\leq i \leq p\,$, the open set  $Z_{i-1}$ is already constructed. Then $Z_{i}$
        is an open neighbourhood of $Y_{i}$ containing $Z_{i-1}$ and such that
        $Z_i\setminus Z_{i-1}$ is  a disjoint union:
        $$
                Z_{i}\setminus Z_{i-1} \cong \coprod_{\si\,i-\mbox{simplex}}\si'\times D_{p-i}\,,
        $$
        where $\si'$ is obtained from $\si$ by a contractive homothety centered at the
        center of $\si$.
\end{itemize}

Let $W_i$ be the inverse image of $Z_i$ by the projection $\V \ar \Tn\,$.
Then $(W_i)_i \,$ is an open filtration of the space $\V$ which satisfies:
$$
        \forall i \in \{1, \cdots ,p\}, W_i\smallsetminus W_{i-1}
        \cong \coprod_{\si\,i-\mbox{simplex}}\O\times \si'\times D_{p-i}.
$$
To the filtration $(W_i)_i$ are associated  two spectral sequences 
converging respectively to $K_*(C(\V))$ and $\hwl$ with $E^2-$tems respectively equal 
to $\oplus_j\cohom_{j}(\zp, C(\O, \Z))$ and $\oplus_j\cohom_{j}(\zp, 
C(\O,\R))\,$. Note that the coefficient in the first homology group comes 
from the identification $K_0(C(\O))\cong C(\O, \Z)$ for totally 
disconnected compact spaces and that the existence of the second 
spectral sequence is a consequence of proposition \ref{longexactsequence}.

\begin{proposition}\label{torsionfree}\cite{ForrestHunton}
The homology group  $\cohom_{*}(\zp, C(\O, \Z))=\oplus_j\cohom_{j}(\zp, C(\O, \Z))$ is 
 torsion free.
 \end{proposition}
    
Recall that there exist a Chern character in Cech cohomology 
$K^*(X)\to \cohom^*(X,\R)$ for every topological space $X\,$. According to remark \ref{cech}, this 
induces a  natural morphism $ch:K^*(U)\to \cohom^*_{\mal,c}(U,\R)$ 
for every open subset $U$ of $\V$ that 
we shall call again the Chern character.

Moreover, according to theorem \ref{intertwin}, the longitudinal 
Chern character induces a morphism between these two spectral  sequences. 
For the $E^2-$term, the longitudinal Chern character $ch_\maL$ and 
the  Chern character $ch$ both induce the  inclusion 
$\oplus_j\cohom_{j}(\zp, C(\O, \Z))\hookrightarrow \oplus_j\cohom_{j}(\zp, 
C(\O,\R))\,$.  
In particular, we get:
\begin{proposition}
    The two morphisms
$ch:K^*(\V)\to \cohom^*_{l}(\V,\R)$ and
$ch_\maL:K^*(\V)\to \cohom^*_{\mal}(\V,\R)$ coincide.
\end{proposition}

    We can check by using spectral sequences associated to the above filtration 
    that  the map $\cohom^*(\V,\R)\to\cohom^*_{\mal}(\V,\R)$ of remark 
    \ref{cech} is an isomorphism.
Actually, we can check that the two above spectral   sequences collapse
at the $E^2$ terms. In particular, since the graduate  associated to 
the filtration $$(\im( K_*(C_0(W_{k}))\to K_*(\Cv)))_{0\leq k\leq p}$$ of 
$K_*(\Cv)$ is 
$\oplus_k\cohom_{k}(\zp, C(\O, \Z))$ and this group being torsion-free by 
proposition \ref{torsionfree}, the abelian group $K_*(C(\V))$ is torsion free.
To be more accurate:

\begin{theorem}\label{homogeneous}\cite{ForrestHunton}
    \begin{itemize}
   \item The morphism $ch:K_*(\V)\to \cohom^*_{\mal}(\V,\R)$ is injective with 
   homogeneous image (i.e its image is graded by the grading of 
   $\cohom^*_{\mal}(\V,\R)\,)$. 
   \item Moreover, the image of $ch:K_*(\V)\to \cohom^*_{\mal}(\V,\R)$ 
   is isomorphic to $\oplus_j\cohom_{j}(\zp,C(\O, \Z))\,$.
   \end{itemize}
   \end{theorem}

To prove the Bellissard conjecture, we shall need the following result of integrality for the top 
dimension Chern character.

\begin{proposition}\label{chp}\cite{BenameurOyono}
Let us denote  by $ ch_{\maL}^k$ the $k-$component of the
longitudinal Chern character
lying in
$\cohom^{k}_{\maL}(\V, \R)$. Then $ \psi_{\Z^p}(ch_{\maL}^p(K_0(\V))$ 
belongs to
$C(\O,\Z)_\zp\,$.
\end{proposition}
\begin{proof}
    Let $x$ be an element of $K_0(\Cv)\,$. According to theorem 
    \ref{homogeneous}, there is an element $x'\in K_0(\Cv)\,$ such 
    that $ch(x)-ch(x')$ is an homogeneous element of degre $p\,$.
    Recall that $Y_{p-1}$ is the $p-1$ skeleton of $\T_p\,$. Let us 
    denote by $V_\O'$ the inverse image of $Y_{p-1}$ under the fibration 
    projection $\V\to\T_p\,$ and by $y$ the image of $x-x'$ under the map 
    $K_0(\Cv)\to K_0(C(V_\O')\,$.
      We can check easily that the Chern
    character  on $V_\O'$ has no component of degre $p$ or higher. In 
    particular, since the Chern character is natural, $ch(y)\in 
    \cohom^*_{\maL}(V_\O',\R) $ vanishes. By using the same trick 
    as we did for $\V\,$, we can check that $K_0(C(V_\O'))$ is 
    torsion-free and that the Chern character of $V_\O'$ is injective.
    In particular, $y=0$ and $x-x'$ lies in the image of
    $K_0(C_0(\V\setminus V_\O'))\to\Kv\,$. We have an identification $\V\setminus 
    V_\O'\cong \O\times]0,1[^p\,$. Since $ch_{\maL}^p(x')=0\,$, under 
    the above identification, $\psi_{\Z^p}(ch_{\maL}^p(x))$ lies in the image of 
    the composition $$K_0(C(\O))\cong K_0(C_0(\O\times ]0,1[^p))\to 
    K_0(\Cv)\stackrel{ch_\maL^p}  {\to}
    \cohom^p_{\maL}(\V,\R)\stackrel{\psi_{\Z^p}}{\cong}C(\O,\R)_\zp\,.$$
    But we have an isomorphism $K_0(C(\O))\cong C(\O,\Z)$ and under this 
    isomorphism, it was shown in \cite{BenameurOyono} that the above 
    composition is equal to the following: $$C(\O,\Z)\to C(\O,\Z)_\zp\to 
    C(\O,\R)_\zp\,.$$
    This conclude the proof.
 \end{proof}

We are now in  position to prove    the Bellissard conjecture.

\begin{theorem}
Let $(\O,\zp)$ be a  dynamical
system with $\O$ a Cantor set and  let $\mu$ be a 
$\zp-$invariant measure on $\O\,$.  Assume that $\O$ has no non-trivial invariant 
compact-open subset. Let $\tau^{\mu}_*$ be the
additive map induced by the trace $\tau^{\mu}$, associated with $\mu$, on the $K$-theory of
the $C^*$-algebra $C(\O) \rtimes \zp$. Then we have:
$$
 \tau^{\mu}_*(K_0(C(\O) \rtimes \zp)) = \Z[\mu].
$$
\end{theorem}

\begin{proof}\
As we said before, we can assume that  $p\,$ is even. From Theorem
\ref{torus}
we deduce that:
$$
 \tau^{\mu}_*(K_0(C(\O) \rtimes \zp)) =
\{\Ind_{\Z^p,\mu}(\de)-\Ind_{\Z^p,\mu}(\pa^{e'}_{\O,\R^p}),
[e] - [e'] \in K_0(C^{\infty,0}(\V))\}.
$$
Using Theorem \ref{measured.index}, we obtain:
$$
 \tau^{\mu}_*(K_0(C(\O) \rtimes \zp)) = <ch^p_{\maL}(K^0(\V)) , [C_{\mu, \Z^p}]>.
$$

Now, from the very definition of $[C_{\mu,\Z^p}]$ and $\psi_{\Z^p}$, we have:
$$
 <x,[C_{\mu, \Z^p}]> = <\psi_{\Z^p}(x),\mu>, \quad \forall x\in
 \cohom^p_{\maL}(\V, \R).
$$
Therefore,
$$
 <ch^p_{\maL}(e),[C_{\mu,\Z^p}]> = <\psi_{\Z^p}  ( ch^p_{\maL}(e)),\mu>.
$$
But according to Proposition \ref{chp}, we have:
$$
        \psi_{\Z^p}(ch^p_{\maL}(e))\in C(\O,\Z)_\zp\,.
$$
Hence:
$$
 \tau^{\mu}_*(K_0(C(\O) \rtimes \zp)) \subset \mu(C(\O,\Z)_\zp)=\Z[\mu].
$$
Since the opposite inclusion is straightforward to check,
 the proof is complete.
\end{proof}

\end{document}